\documentclass[reqno]{amsart}

\usepackage{amsmath}
\usepackage{amssymb}
\usepackage[all]{xy}
\usepackage{fullpage}
\numberwithin{equation}{section}
\newcommand{\bC}{\mathbb{C}}

\newcommand{\bP}{\mathbb{P}}
\newcommand{\bQ}{\mathbb{Q}}

\newcommand{\bZ}{\mathbb{Z}}

\newcommand{\cA}{\mathcal{A}}

\newcommand{\cD}{\mathcal{D}}

\newcommand{\wt}{\widetilde}

\newcommand{\nd}{\noindent}
\newcommand{\iso}{\simeq}
\renewcommand{\proof}{\nd \textbf{Proof}.\,\,}
\newcommand{\Spec}{\text{Spec}\,}

\newcommand{\Pic}{\text{Pic}}

\newcommand{\D}{\Delta}

\makeatletter
\renewcommand{\section}{\@startsection{section}{1}{0mm}{1.5\baselineskip}{0.5\baselineskip}{\Large\bf\center}}
\renewcommand{\subsection}{\@startsection{subsection}{2}{0mm}{0.5\baselineskip}{-0.5em}{\bf}}
\renewcommand{\subsubsection}{\@startsection{subsubsection}{3}{0mm}{0.5\baselineskip}{-0.5em}{\bf}}

\makeatother

\begin{document}

\begin{center}
{\LARGE \bf  Dubrovin's conjecture for $IG(2,6)$}

\vspace{10pt}

{\large Sergey~Galkin, Anton~Mellit, Maxim~Smirnov}

\vspace{2pt}

\end{center}

\vspace{10pt}

\textsc{Abstract.} We show that the big quantum cohomology of the symplectic isotropic Grassmanian $IG(2,6)$ is generically semisimple, whereas its small quantum cohomology is known to be non-semisimple. This gives yet another case where Dubrovin's conjecture holds and stresses the need to consider the big quantum cohomology in its formulation.

\section{Introduction}

The main purpose of this paper is to give an explicit example of a smooth projective variety $X$ such that its small quantum cohomology $qH^*(X)$ is not generically semisimple, whereas its big quantum cohomology $QH^*(X)$ is generically semisimple. Namely, we show that this pattern holds for the symplectic isotropic Grassmanian $IG(2,6)$. 

In general the small quantum cohomology of isotropic Grassmanians is studied in \cite{BuKrTa}, and an explicit presentation in terms of generators and relations is given. Based on this presentation in \cite[Sec.~7]{ChMaPe} it is shown that $qH^*(IG(2,6))$ is not generically semisimple.

The Grassmanian $IG(2,6)$ appears to be the simplest explicit example of this sort available in the literature. More examples will appear in a forthcoming paper by Nicolas Perrin \cite{Pe}, who also independently proved the generic semisimplicity of $QH^*(IG(2,6))$ simultaneously with us.

\subsection{Remark} Throughout the text we use notation $QH^*(X)$ for the big quantum cohomology of $X$ and $qH^*(X)$ for the small quantum cohomology of $X$. We should mention that this notation is not standard, but it is very convenient for our purposes. To avoid any confusion, in Section~\ref{Sec.: Reminder on QH} we will briefly recall both concepts, since they are crucial for this paper.

\subsection{Dubrovin's conjecture}

Dubrovin's conjecture (see \cite{Du}) gives an intriguing relation between the quantum cohomology of a smooth projective variety $X$ and the bounded derived category of coherent sheaves on it. Namely, it says that the generic semisimplicity of $QH^*(X)$ is equivalent to the existence of a full exceptional collection in $\cD^b(X)$. Here $\cD^b(X)$ denotes the bounded derived category of coherent sheaves on~$X$.

The main motivation for studying examples like $IG(2,6)$ comes from the conjecture above. In \cite{Ku} it is shown that $\cD^b(IG(2,6))$ has a full exceptional collection. Hence, it is expected that $QH^*(IG(2,6))$ is generically semisimple, whereas $qH^*(IG(2,6))$ is not. 

For completeness let us just mention that one expects to find many more such examples among homogeneous spaces $G/P$, where $G$ is a semi-simple algebraic group and $P$ is a parabolic subgroup. The reason is that $\cD^b(G/P)$ is conjectured to have a full exceptional collection (see \cite[Conj.~1.1]{KuPo}). In particular, for arbitrary symplectic isotropic Grassmanians $IG(k,2n)$ the big quantum cohomology is expected to be generically semisimple, and one can try to find a pattern when $qH^*(IG(k,2n))$ is not. For results in this direction related to semisimplicity we refer to \cite[Sec.~6]{ChPe}.

We should remark that the above formulation of Dubrovin's conjecture is not complete. Here we consider its first part only. The remaining two parts (see \cite[Conj.~4.2.2]{Du} and \cite{GaGoIr}) are not considered in this paper. Let us just mention that for Grassmanians $G(k,n)$ all three parts are known to hold (see \cite{GaGoIr, Ue} and references therein). For a general introduction to this topic we refer to \cite{Du, Ba, GaGoIr, HeMaTe}.

\subsubsection{Some known examples}

Here we list some instances where Dubrovin's conjecture is known to hold. The simplest example is provided by the projective spaces $\bP^n$. Indeed, it is well-known that $qH^*(\bP^n)$ is generically semisimple (for example, see \cite{KoMa}), and $\cD^b(\bP^n)$ has a full exceptional collection (see \cite{Be}). 

One way to generalize the example of $\bP^n$ is to look at Grassmanians $G(k,n)$. By \cite[Prop.~6.5]{Ab} (or \cite[Cor.~8]{GaGo}) it is known that $qH^*(G(k,n))$  is generically semisimple, and by \cite{K} we know that $\cD^b(G(k,n))$ has a full exceptional collection. 

Another way to generalize the example of $\bP^n$ is to consider arbitrary toric varieties. For a smooth projective toric variety $X$ (no Fano assumption!) it is known that $QH^*(X)$ is generically semisimple (see \cite{Ir}). Note that from loc.cit. it is not known whether $qH^*(X)$ is generically semisimple for an arbitrary toric variety $X$, but only when $X$ is a weak Fano variety. A full exceptional collection in $\cD^b(X)$ exists for any toric $X$ by \cite{Ka}.

In the above examples already the small quantum cohomology is generically semisimple. As far as we know $IG(2,6)$ is the first explicit example in the literature where one really needs the big quantum cohomology in order to formulate Dubrovin's conjecture.

\medskip

\textbf{Acknowledgements.} We are very grateful to Tarig Abdelgadir, Boris Dubrovin, Vasily Golyshev, Vassily Gorbounov, Hiroshi Iritani, Alexander Kuznetsov, Yuri Manin, Nicolas Perrin, and Evgeny Shinder for valuable discussions, comments, and attention to this work.

\section{Reminder on quantum cohomology}
\label{Sec.: Reminder on QH}

Here we briefly recall the definition of the quantum cohomology ring for a smooth projective variety $X$. It is convenient for us to define it in a fixed basis of $H^*(X, \bQ)$ but, of course, the resulting ring structure is independent of it. We assume some familiarity with the subject and do not go into much detail. For a thorough introduction we refer to \cite{Ma}. 

\subsection{Definition}

Let $X$ be a smooth projective variety and assume for simplicity that $H^{odd}(X,\bQ)=0$. Fix a graded basis $\Delta_0, \dots , \Delta_n$ in $H^*(X, \bQ)$ and dual linear coordinates $t_0, \dots, t_n$. It is customary to choose $\D_0=1$. Let $N_X$ be the Novikov ring of $X$, i.e. the ring of formal power series in $q^{\beta}$, where $\beta$ runs over the Mori cone of curves on $X$.

The genus zero Gromov-Witten potential of $X$ is an element $\Phi \in N_X[[t_0, \dots, t_n]]$, i.e. a formal power series in variables $t_0, \dots, t_n$ with coefficients in the Novikov ring $N_X$, defined by the formula
\begin{align}\label{Eq.: GW potential}
&\Phi = \sum_{(i_0, \dots , i_n)}  \langle \Delta_0^{\otimes i_0}, \dots, \Delta_n^{\otimes i_n} \rangle \frac{t_0^{i_0} \dots t_n^{i_n} }{i_0!\dots i_n!}, 
\end{align}
where $\langle \Delta_0^{\otimes i_0}, \dots, \Delta_n^{\otimes i_n} \rangle := \sum_{\beta} \langle \Delta_0^{\otimes i_0}, \dots, \Delta_n^{\otimes i_n} \rangle_{\beta} q^{\beta}$, and $\beta$ runs over the cone of effective curve classes. The rational numbers $\langle \Delta_0^{\otimes i_0}, \dots, \Delta_n^{\otimes i_n} \rangle_{\beta}$ are Gromov-Witten invariants of $X$ of the curve class $\beta$.

Using \eqref{Eq.: GW potential} one defines the quantum cohomology ring of $X$. Namely, on the basis elements we put
\begin{align}\label{Eq.: Quantum cohomology}
\Delta_a \star \Delta_b = \sum_c \Phi_{abc} \Delta^c,
\end{align}
where $\Phi_{abc}=\frac{\partial^3 \Phi}{\partial t_a \partial t_b \partial t_c}$, and $\D^0, \dots, \D^n$ is the basis dual to  $\D_0, \dots, \D_n$ with respect to the Poincar\'e pairing. Expression \eqref{Eq.: Quantum cohomology} is naturally interpreted as an element of $H^*(X, \bQ)\otimes_{\bQ} N_X[[t_0, \dots, t_n]]$.

\smallskip

It is well known that \eqref{Eq.: Quantum cohomology} makes $H^*(X, \bQ)\otimes_{\bQ} N_X[[t_0, \dots, t_n]]$ into a commutative, associative, graded $N_X[[t_0, \dots, t_n]]$-algebra with the identity element $\D_0$. We will denote this algebra $QH^*(X)$. For convenience we recall the definition of the grading:
\begin{align*}
& \deg(\D_i)=\frac{|\D_i|}{2}, \quad \deg(q^{\beta})=(-K_X,\beta), \quad \deg(t_i)=1-\frac{|\D_i|}{2},
\end{align*}
where $|\D_i|$ is the cohomological degree of $\D_i$. 

Sometimes the algebra $QH^*(X)$ is called the \textit{big quantum cohomology algebra} of $X$ to distinguish it from a simpler object called the \textit{small quantum cohomology algebra}. It is the quotient of $QH^*(X)$ with respect to the ideal $(t_0, \dots, t_n)$. We will denote the latter $qH^*(X)$ and use $\circ$ instead of $\star$ for the product in this algebra. 

It is equivalent to say that $qH^*(X)=H^*(X,\bQ) \otimes_{\bQ} N_X$ as a vector space, and the $N_X$-algebra structure is defined by putting $\Delta_a \circ \Delta_b = \sum_c \langle \D_a, \D_b, \D_c \rangle \Delta^c$.

\smallskip

\subsubsection{Remark} \textit{(i)} Often in the literature one uses $QH^*(X)$ both for big and small quantum cohomology. To avoid any confusion we have decided to use the non-standard notation described above. 

\textit{(ii)} At first sight our definition of the small quantum cohomology, i.e. the ring structure defined by 3-pointed GW invariants, differs slightly from \cite{Ma}. It is not hard to see that both structures are equivalent. 

\textit{(iii)} In Section \ref{Sec.: Deformation and semisimplicity} we will use $\star$ to denote the big quantum product in a more restricted sense. Namely, instead of a full deformation we will consider only one-parameter deformation.

\textit{(iv)} Throughout the text we work only with the even cohomology. Therefore, all degrees appearing below are Chow-ring degrees.

\subsection{Lemma}
\label{SubSec.: Multiplication by h}

Let $\Pic X=\bZ$, $h \in H^2(X,\bQ)$ and $\gamma \in H^*(X,\bQ)$. Then
\begin{align*}
\frac{\partial}{\partial t_i} ( h \star \gamma ) = (h,\beta) \, q\frac{d}{dq} (\Delta_i \star \gamma), 
\end{align*}
where $q=q^{\beta}$ with $\beta$ being the generator of the cone of effective curves.

\smallskip

\proof We have the chain of equalities
\begin{align*}
\frac{\partial}{\partial t_i} ( h \star \gamma ) = \nabla_{\Delta_i}( h \star \gamma ) = \nabla_h (\Delta_i \star \gamma) = (h,\beta) \, q\frac{d}{dq} (\Delta_i \star \gamma),  
\end{align*}
where the first one holds by definition, the second follows from potentiality of the quantum product (see the proof of Proposition 1.6 in \cite{Ma}), and the last one holds due to the divisor axiom for GW invariants.$\blacksquare$

\subsection{Corollary}

For an element $\alpha$ in  $H^*(X, \bQ)$ let $M_{\alpha}$ and $\wt{M}_{\alpha}$ be the matrices of the small and the big quantum multiplication by $\alpha$ in the basis $\D_0, \dots, \D_n$. Moreover, to simplify notation, put $M_i=M_{\D_i}$ and $\wt{M}_i=\wt{M}_{\D_i}$.

From Lemma~\ref{SubSec.: Multiplication by h} it follows that for $h \in H^2(X,\bQ)$ we have
\begin{align}\label{Eq.: Integral}
\wt{M}_h = M_h + (h, \beta) q\frac{d}{dq} \sum_i \int_0^{t_i} \wt{M}_i \,\, dt_i,
\end{align} 
and restricting to the first infinitesimal neighbourhood we get
\begin{align}\label{Eq.: Multiplication by h in the first nbhd}
\wt{M}_h = M_h + (h, \beta) q \frac{d}{dq} \sum_i t_i M_i + O(t_at_b).
\end{align}

\section{Small quantum cohomology of $IG(2,6)$}

The aim of this paragraph is to fix notation and to describe the structure of the small quantum cohomology ring of $IG(2,6)$ following \cite{ChMaPe}.

\subsection{Geometry of $IG(2,6)$} 

Let $V$ be a six-dimensional complex vector space and $\omega \in \Lambda^2 V^*$ a symplectic form. The Grassmanian $IG(2,V)$ parametrizes two-dimensional isotropic subspaces of $V$. It is a smooth projective variety of dimension $\dim X=7$.

One can give an explicit presentation of $IG(2,V)$ as a hyperplane section of $G(2,V)$. Indeed, consider the Grassmanian $G(2,V)$, and its Pl\"ucker embedding into $\bP(\Lambda^2 V) \iso \bP^{14}$. The form $\omega$ defines a hyperplane in $\bP(\Lambda^2 V)$. The section of $G(2,V)$ by this hyperplane is $X=IG(2,6)$. 

Similarly to the case of usual Grassmanians, there exists a basis of cohomology indexed by the so called 1-strict partitions $\lambda$ (see \cite{ChMaPe} and references therein). The cohomology class corresponding to such a partition $\lambda$ will be denoted $\Delta_{\lambda}$. Explicitly, there are 12 such classes
\begin{align}\label{Eq.: Basis of cohomology}
\D_0; \quad \D_1; \quad \D_2, \D_{1,1}; \quad \D_3, \D_{2,1}; \quad \D_4, \D_{3,1}; \quad \D_{4,1}, \D_{3,2}; \quad \D_{4,2}; \quad \D_{4,3},
\end{align}
and the Chow-ring degree of $\Delta_{\lambda}$ is $|\lambda|=\lambda_1+\lambda_2$.

\subsection{Gromov-Witten invariants}

Recall that $\Pic \, X=\bZ$ and let $l$ be the generator of the cone of effective curves on $X$. Since $-K_X=5\D_1$, the dimension axiom for an n-pointed GW invariant $\langle \gamma_1 , \dots , \gamma_n \rangle$ becomes
\begin{align*}
n+4 + 5d = \sum_i |\gamma_i|,
\end{align*} 
and implies that $\langle \D_2 , \D_2 , \D_2 \rangle= \langle \D_2 , \D_2 , \D_2, \D_2 \rangle = 0$.

\subsection{Multiplication table}

As in \cite[Sec.7]{ChMaPe} we represent the small quantum cohomology of $X$ specialized at $q=1$ by the table
\begin{align}\label{Eq.: Table of Small QH}
\begin{array}{cccccccccccc}
& \D_1 & \D_2 & \D_{1,1} & \D_{2,1} & \D_3 & \D_{3,1} & \D_4 & \D_{4,1} & \D_{3,2} & \D_{4,2} &\D_{4,3} \\
& & & & & & & & & & &\\
Z_0 & 0 &\varepsilon & -\varepsilon & 0 & 0 & 0 & 0 & -1 & 1 & 0 & -\varepsilon       \\
Z_1 & s & 0 & s^2 & s^3 & -s^3 & s^4 & -s^4 & 0 & -1 & -s         & -s^2           \\
Z_2 & t & \frac{2}{3}t^2 & \frac{1}{3}t^2 & \frac{1}{3}t^3& \frac{1}{3}t^3& \frac{1}{9}t^4 & \frac{1}{9}t^4 & 2 & 1 & t & \frac{1}{3}t^2  \\
\end{array}
\end{align}

Some explanations are in order. The small quantum cohomology of $X$ specialized at $q=1$ is a finite dimensional algebra $A$. Explicitly it is given as a direct product 
\begin{align}\label{Eq.: Ring of Small QH}
&A=A_0 \times A_1 \times A_2,\\ \notag
&A_0 = \bC[\varepsilon]/ \varepsilon^2 , \quad A_1=\bC[s]/(s^5+1) , \quad A_2=\bC[t]/(t^5-27).
\end{align}
 Therefore, its spectrum $\Spec(A)$ is a disjoint union of three components $Z_i=\Spec(A_i)$. Each element $\D_{\lambda}$ is a function on $Z=\Spec(A)$, and we can look at its values at different points. Table \eqref{Eq.: Table of Small QH} collects these values grouped according to components $Z_i$.

From \eqref{Eq.: Table of Small QH} one can reconstruct all small quantum products. Moreover, using the fact that the quantum multiplication is graded, one can also put the parameter $q$ back. For example, the small quantum product with $\D_1$ is 
\begin{align}\label{Eq.: Small multiplication table with Delta_1}
\begin{array}{ll}
& \D_1 \circ \D_1 = \D_2+\D_{1,1} \\
& \D_1 \circ \D_2 = \D_3+\D_{2,1} \\
& \D_1 \circ \D_{1,1} = \D_{2,1} \\
& \D_1 \circ \D_3 = 2 \D_4 + \D_{3,1}\\
& \D_1 \circ \D_{2,1} =  \D_4 + 2\D_{3,1}\\
& \D_1 \circ \D_4 = \D_{4,1}+q\D_0 \\
& \D_1 \circ \D_{3,1} = \D_{4,1}+\D_{3,2} \\
& \D_1 \circ \D_{4,1} = \D_{4,2}+q\D_1 \\
& \D_1 \circ \D_{3,2} = \D_{4,2} \\
& \D_1 \circ \D_{4,2} = \D_{4,3}+q\D_2 \\
& \D_1 \circ \D_{4,3} = q\D_3.
\end{array}
\end{align}
Analogously one obtains multiplication tables with $\D_2$ and $\D_3$
\begin{align}\label{Eq.: Small multiplication table with Delta_2 and Delta_3}
\begin{array}{lll}
\D_2 \circ \D_2 = 2 (\D_4 + \D_{3,1})                    &\quad\quad  & \D_3 \circ \D_2 = 2\D_{4,1}+\D_{3,2}+q \D_0      \\
\D_2 \circ \D_{1,1} = \D_4 + \D_{3,1}                    &\quad\quad  & \D_3 \circ \D_{1,1} = \D_{4,1} +q \D_0           \\
\D_2 \circ \D_3 = 2 \D_{4,1} + \D_{3,2} + q \D_0         &\quad\quad  & \D_3 \circ \D_3 = 2\D_{4,2} + q\D_1              \\
\D_2 \circ \D_{2,1} = 2 \D_{4,1} + \D_{3,2} + q \D_0     &\quad\quad  & \D_3 \circ \D_{2,1} = \D_{4,2} + 2q\D_1          \\
\D_2 \circ \D_4 = \D_{4,2} + q \D_1                      &\quad\quad  & \D_3 \circ \D_4 = \D_{4,3} + q \D_2              \\
\D_2 \circ \D_{3,1} = \D_{4,2} + q \D_1                  &\quad\quad  & \D_3 \circ \D_{3,1} = q (\D_2 + \D_{1,1})        \\
\D_2 \circ \D_{4,1} = \D_{4,3} + q (\D_2 +\D_{1,1})      &\quad\quad  & \D_3 \circ \D_{4,1} = q(\D_{2,1}+\D_3)           \\
\D_2 \circ \D_{3,2} = q \D_2                             &\quad\quad  & \D_3 \circ \D_{3,2} = q \D_{2,1}                 \\
\D_2 \circ \D_{4,2} = q (\D_3 + \D_{2,1})                &\quad\quad  & \D_3 \circ \D_{4,2} = q (2\D_{3,1} + \D_4)       \\
\D_2 \circ \D_{4,3} = q (\D_4 + \D_{3,1})                &\quad\quad  & \D_3 \circ \D_{4,3} = q (\D_{4,1} + \D_{3,2}).
\end{array}
\end{align}

Formulas \eqref{Eq.: Small multiplication table with Delta_1} and \eqref{Eq.: Small multiplication table with Delta_2 and Delta_3} suffice to reconstruct the whole multiplication table. Indeed, if we let $M_{\lambda}$ be the matrix of the small quantum multiplication with $\D_{\lambda}$, then 
\begin{align*}
\begin{array}{ll}
M_{1,1}=M_1^2-M_2 \\
M_{2,1}=M_1M_2-M_3 \\
M_{3,1}=-\frac{1}{3}M_1(M_3-2M_{2,1}) \\
M_4=M_1M_{2,1}-2M_{3,1} \\
M_{4,1}=M_1M_4-qM_0 \\
M_{3,2}=M_1M_{3,1}-M_{4,1} \\
M_{4,2}=M_1M_{4,1}-qM_1 \\
M_{4,3}=M_1M_{4,2}-qM_2. 
\end{array}
\end{align*}

\subsection{Nilpotent}

It is clear from \eqref{Eq.: Ring of Small QH} that $qH^*(X)$ specialized at $q=1$ is not semisimple. Indeed, $Z$ has 10 reduced points and one fat point. The corresponding nilpotent (see \cite{ChMaPe}), which is unique up to a constant factor, is 
\begin{align}\label{Eq.: Nilpotent specialized at q=1}
\D_{4,3} -  \D_2 +  \D_{1,1}.
\end{align}
Moreover, it is not hard to see that \eqref{Eq.: Nilpotent specialized at q=1} is a specialization of a nilpotent defined globally over $q$. Indeed, it is given by $c_0=\D_{4,3} - q \D_2 + q \D_{1,1}$.

\section{Deformation along $\D_2$ and semisimplicity}
\label{Sec.: Deformation and semisimplicity}

\subsection{Theorem}

The big quantum cohomology of $IG(2,6)$ is generically semisimple.

\subsection{Proof}

As we will see, to exhibit the generic semisimplicity it is enough to consider the deformation of the small quantum cohomology along $\D_2$ instead of the full big quantum cohomology $QH^*(X)$. Let $\cA$ be this deformation, i.e. $\cA$ is the quotient of $QH^*(X)$ with respect to the ideal generated by coordinates $t_{\lambda}$ for $\lambda \neq 2$. To simplify the notation we will write $t$ instead of $t_2$. The multiplication in $\cA$ will be denoted $\star$.

To show the generic semisimplicity of $\cA$ (and, hence, of $QH^*(X)$) we exhibit an element that has distinct eigenvalues. Namely, we show that the element $\gamma=\D_2 + \D_1$ satisfies this property. The strategy can be described as follows. First, using Lemma~\ref{SubSubSec.: Multiplication in the 1st nbhd} below, we compute the matrix $\wt{M}_{\gamma}$ of multiplication by $\D_2 + \D_1$ in $\cA$ modulo terms quadratic in $t$. Secondly, we compute the characteristic polynomial $P(x)$ of $\wt{M}_{\gamma}$. Finally, by looking at Newton polygons of $P(x)$ and $P_x'(x)$ one deduces that the roots of $P(x)$ are all distinct.

The computations outlined above are done using PARI/GP and the source code is available at \cite{Source}. In the final step we use a well-known relation between roots of a polynomial over a discrete valuation ring (in our case it is the ring of formal power series in $t$) and its Newton polygon (for example, see \cite[IV.3]{Ko}). Namely, in PARI/GP we obtain 
\begin{align}\label{Eq.: Charpoly(HH+DD)}
&P(x)=P_0(x)+tP_1(x)+O(t^2),
\end{align}
where
\begin{align*}
&P_0(x)=x^{12} - 60qx^9 - 90qx^8 -(96q^2 + 26q)x^7 - 60q^2x^4 - 90q^2x^3 -(96q^3 + 27q^2)x^2\\ 
&P_1(x)=-30qx^{10} -  96qx^9 - 36qx^8 + 152q^2x^7 + 120q^2x^6 + (-32q^3 + 186q^2)x^5 + \\ 
& \hspace{35pt} + (240q^3 - 26q^2)x^4  - 36q^2x^3 + 152q^3x^2 - 30q^3x -32q^4 - 9q^3.
\end{align*}

The characteristic polynomial of $\D_2 + \D_1$ in the small quantum cohomology is equal to the reduction of \eqref{Eq.: Charpoly(HH+DD)} modulo $t$, i.e. to $P_0(x)$. One can check that $P_0(x)$ has ten simple roots and one root of multiplicity two.

The Newton polygon of \eqref{Eq.: Charpoly(HH+DD)} is
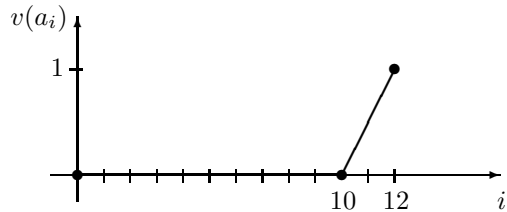
$$
\begin{picture}(100,85)
\put(-20,20){\vector(1,0){170}}
\put(-10,10){\vector(0,1){70}}
\put(-10,20){\circle*{4}}
\put(90,20){\circle*{4}}
\put(110,60){\circle*{4}}

\put(90,10){\makebox(1,1){$10$}}
\put(110,10){\makebox(1,1){$12$}}
\put(150,10){\makebox(1,1){$i$}}
\put(-25,79){\makebox(1,1){$v(a_i)$}}
\put(-18,60){\makebox(1,1){$1$}}

\put(-13,60){\line(1,0){5}}

\put(0,17){\line(0,1){5}}
\put(10,17){\line(0,1){5}}
\put(20,17){\line(0,1){5}}
\put(30,17){\line(0,1){5}}
\put(40,17){\line(0,1){5}}
\put(50,17){\line(0,1){5}}
\put(60,17){\line(0,1){5}}
\put(70,17){\line(0,1){5}}
\put(80,17){\line(0,1){5}}
\put(100,17){\line(0,1){5}}
\put(110,17){\line(0,1){5}}

\thicklines
\put(90,20){\line(1,2){20}}
\put(-10,20){\line(1,0){100}}
\end{picture}
$$
where $P(x)=a_0x^{12}+a_1x^{11}+\dots + a_{11}x+a_{12}$ and $v$ is the valuation with respect to $t$. Using Lemma 4 of \cite[IV.3]{Ko} we immediately obtain that $P(x)$ has ten roots of valuation zero, i.e. their expansions start from a constant term, and two roots of valuation $\frac{1}{2}$. Since $P_0(x)$ has ten simple roots, ten roots of $P(x)$ of valuation zero are distinct. Considering the Newton polygon for $P'_x(x)$ we see that it has ten roots of valuation zero and one root of valuation one. This implies that $P(x)$ and $P'_x(x)$ have no common roots. Therefore, all roots of $P(x)$ are distinct.

\subsubsection{Lemma}
\label{SubSubSec.: Multiplication in the 1st nbhd}

Multiplication in $\cA/t^2$ can be reconstructed from $qH^*(X)$, i.e. from 3-pointed Gromov-Witten invariants.

\smallskip

\proof Consider elements 
\begin{align}\label{Eq.: Cool basis}
& f_i=h^{\star i} \quad \text{for $0 \leq i \leq 10$}   \\ \notag
& f_{11}=\D_2.
\end{align}
These elements form a basis of $\cA$ because their reductions modulo $t$ are equal to $1, h, h^{\circ 2}, \dots, h^{\circ 10}, \D_2$, which form a basis in $\cA/t = qH^*(X)$.

In Basis \eqref{Eq.: Cool basis} the multiplication table looks very simple. Clearly, if $i+j \leq 10$, we have $f_i \star f_j = f_{i+j}$. Therefore, the only non-trivial part are the products $f_i \star f_{11}$.

For $i\neq 11$ we can compute the products $f_i \star f_{11}$ in the first neighbourhood due to Lemma~\ref{SubSec.: Multiplication by h}. So the only possibly non-trivial product is $f_{11} \star f_{11}=\D_2 \star \D_2$. 

To compute $\D_2 \star \D_2$ we do the following. First we write it as
\begin{align*}
&\D_2 \star \D_2 = \sum_{i=0}^{10} (\D_2 \star \D_2, f_i)f^i + (\D_2 \star \D_2, f_{11})f^{11},
\end{align*}
where $f^i$ is the dual of $f_i$ with respect to the Poincar\'e pairing. Further, note that 
\begin{align*}
&(\D_2 \star \D_2, f_i)=(\D_2 , \D_2 \star f_i) =(\D_2 \star f_i, \D_2),
\end{align*}
so $(\D_2 \star \D_2, f_i)$ can be extracted from the products $f_i \star \D_2$, which we already know for $i \neq 11$. Finally, for the remaining $(\D_2 \star \D_2, f_{11})$ we have 
\begin{align*}
&(\D_2 \star \D_2, f_{11})=(\D_2 \star \D_2, \D_2)=\langle \D_2 , \D_2 , \D_2 \rangle + \langle \D_2 , \D_2 , \D_2, \D_2 \rangle t +O(t^2).
\end{align*}
By the dimension axiom for Gromov-Witten invariants we have $\langle \D_2 , \D_2 , \D_2 \rangle= \langle \D_2 , \D_2 , \D_2, \D_2 \rangle = 0$. Therefore, we do not need to compute any 4-pointed Gromov-Witten invariants to get the multiplication table in the first neighbourhood along $t$.  $\blacksquare$

\subsection{Spectrum of the Euler field}
\label{SubSec.: Spectrum of Euler}

Here we will show that the operator of multiplication with the standard Euler vector field $E$ of $QH^*(X)$ has simple spectrum (see \cite[Proposition 5.3.4]{Ma} for a definition). It is enough to show this for the restriction of $E$ to the $\D_2$-direction, i.e. to consider the element
\begin{align}\label{Eq.: Euler field}
E=5 \Delta_1 -t \Delta_2,
\end{align}
where used $E$ to denote this restriction as well.

Combining \eqref{Eq.: Integral} and Lemma~\ref{SubSubSec.: Multiplication in the 1st nbhd} one can compute the matrix $\wt{M}_E$ in the second neighbourhood along $t$. Indeed, from Lemma~\ref{SubSubSec.: Multiplication in the 1st nbhd} we know $\wt{M}_2$ modulo $t^2$. Therefore, we know $t\wt{M}_2$ modulo $t^3$. Moreover, by \eqref{Eq.: Integral} 
\begin{align*}
\wt{M}_1 = M_1 + q\frac{d}{dq} \int_0^{t} \wt{M}_2 \,\, dt.
\end{align*}
Since we know $\wt{M}_2$ modulo $t^2$, this formula gives us $\wt{M}_1$ modulo $t^3$. Hence, we obtain $\wt{M}_E$ modulo $t^3$.

Implementing the above algorithm in PARI/GP and computing the characteristic polynomial $P(x)$ of $\wt{M}_E$ we obtain\footnote{Note that here the precision of the characteristic polynomial is higher than that of $\wt{M}_E$.}
\begin{align}\label{Eq.: Charpoly(E)} 
&P(x)=P_0(x)+tP_1(x)+t^2P_2(x) - 39062500q^3t^3 + xO(t^3) + O(t^4),
\end{align}
where
\begin{align*}
&P_0(x)=x^{12} - 81250qx^7 - 263671875q^2x^2 \\
&P_1(x)= -11250qx^8 - 35156250q^2x^3 \\
&P_2(x)= -900q x^9 - 78125q^2x^4.
\end{align*}
Arguing the same way as after \eqref{Eq.: Charpoly(HH+DD)} gives that \eqref{Eq.: Charpoly(E)} has distinct roots. Hence, the standard Euler vector field has simple spectrum.

\end{document}